\documentclass{article}

\usepackage{amssymb}
\usepackage[dvips]{graphicx}

\newtheorem{theorem}{Theorem}

\begin{document}

\title{Two multicolor Ramsey numbers}

\author{Alexander Engstr\"om\footnote{
Institute of Theoretical Computer Science, ETH Z\"urich, CH-8092 Z\"urich, Switzerland \break
engstroa@inf.ethz.ch \break
Research supported by ETH and Swiss National Science Foundation Grant PP002-102738/1}
}

\maketitle

\begin{quote}
Two new bounds for multicolor Ramsey numbers are proved:\break 
\mbox{$R(K_3,K_3,C_4,C_4)\geq 27$} and $R_4(C_4)\leq 19$.
\end{quote}

\section{Introduction}
We prove two new bounds for multicolor Ramsey numbers, a lower
bound for $R(K_3,K_3,C_4,C_4)$ by coloring $K_{26}$, and an upper
bound for $R_4(C_4)$ by a density argument.

\section{The Ramsey number $R(K_3,K_3,C_4,C_4)$}
From the survey of Ramsey numbers by Radziszowski \cite{r} we know
that\break $R(K_3,K_3,C_4,C_4)\geq 26$. We use 
$C_5$-decompositions to construct a four-coloring of the edges of
$K_{26}$, which show that $R(K_3,K_3,C_4,C_4) \geq 27$. The
technique used in this section was invented by Exoo and Reynolds
\cite{er}.

\begin{theorem}
$R(K_3,K_3,C_4,C_4) \geq 27$.
\end{theorem}

\noindent
\textbf{Proof:} 
Let $X, Y, I, \bar{0}, $ and $\bar{1}$ be defined by

\begin{center}
\begin{tabular}
{p{3pt}p{55pt}p{3pt}p{55pt}p{3pt}p{55pt}p{3pt}p{15pt}p{3pt}p{15pt}}
\emph{X}=&$\left[ 
\begin{array}{p{0pt}p{0pt}p{0pt}p{0pt}p{0pt}}
0&1&0&0&1\\
1&0&1&0&0\\
0&1&0&1&0\\
0&0&1&0&1\\
1&0&0&1&0\\
\end{array}
\right]$
&
\emph{Y}=&$\left[
\begin{array}{p{0pt}p{0pt}p{0pt}p{0pt}p{0pt}}
0&0&1&1&0\\
0&0&0&1&1\\
1&0&0&0&1\\
1&1&0&0&0\\
0&1&1&0&0\\
\end{array}
\right]$
& 
\emph{I}=&$\left[
\begin{array}{p{0pt}p{0pt}p{0pt}p{0pt}p{0pt}}
1&0&0&0&0\\
0&1&0&0&0\\
0&0&1&0&0\\
0&0&0&1&0\\
0&0&0&0&1\\
\end{array}
\right]$
&
$\bar{0}$=&$\left[
\begin{array}{p{0pt}}
0\\
0\\
0\\
0\\
0\\
\end{array}
\right]$
&
$\bar{1}$=&$\left[
\begin{array}{p{0pt}}
1\\
1\\
1\\
1\\
1\\
\end{array}
\right]$
\\
\end{tabular}
\end{center}

The critical coloring which shows that $R(K_3,K_3)>5$ and 
$R(C_4,C_4)>5$ have the adjacency matrices $X$ and $Y$. 
Observe that $X+Y+I$ is the all-ones $5\times 5$ matrix. 

We now construct four $26\times 26$ adjacency matrices $M_i$, so
that $M_1$ and $M_2$ contain no $K_3$, and $M_3$ and $M_4$ contain
no $C_4$.

Given a triangle-free graph on $n$ vertices, we can construct a
triangle-free graph on $nm$ vertices by replacing each vertex with
$m$ vertices and each edge with $K_{m,m}$. We construct the two
first graphs, which are isomorphic, by beginning with $C_5$,
replacing the edges with $K_{5,5}-e$, and then adding a vertex
with five edges.

\[
\begin{array}{cc}
M_1 = \left[ \begin{array}{cccccc}
0 & X & X & X & X & \bar{1} \\
X & X & X & X & X & \bar{0} \\
X & X & X & X & X & \bar{0} \\
X & X & X & X & X & \bar{0} \\
X & X & X & X & X & \bar{0} \\
\bar{1}^T  & \bar{0}^T & \bar{0}^T  & \bar{0}^T & \bar{0}^T & 0
\end{array} \right]
&
M_2 = \left[ \begin{array}{cccccc}
Y & Y & Y & Y & Y & \bar{0} \\
Y & Y & Y & Y & Y & \bar{0} \\
Y & Y & Y & Y & Y & \bar{0} \\
Y & Y & Y & Y & Y & \bar{0} \\
Y & Y & Y & Y & 0 & \bar{1} \\
\bar{0}^T  & \bar{0}^T & \bar{0}^T  & \bar{0}^T & \bar{1}^T & 0
\end{array} \right]
\end{array}
\]
\begin{figure}[!ht]
   \begin{center}
     \begin{tabular}{cc}
      \includegraphics*[width=0.407\textwidth]{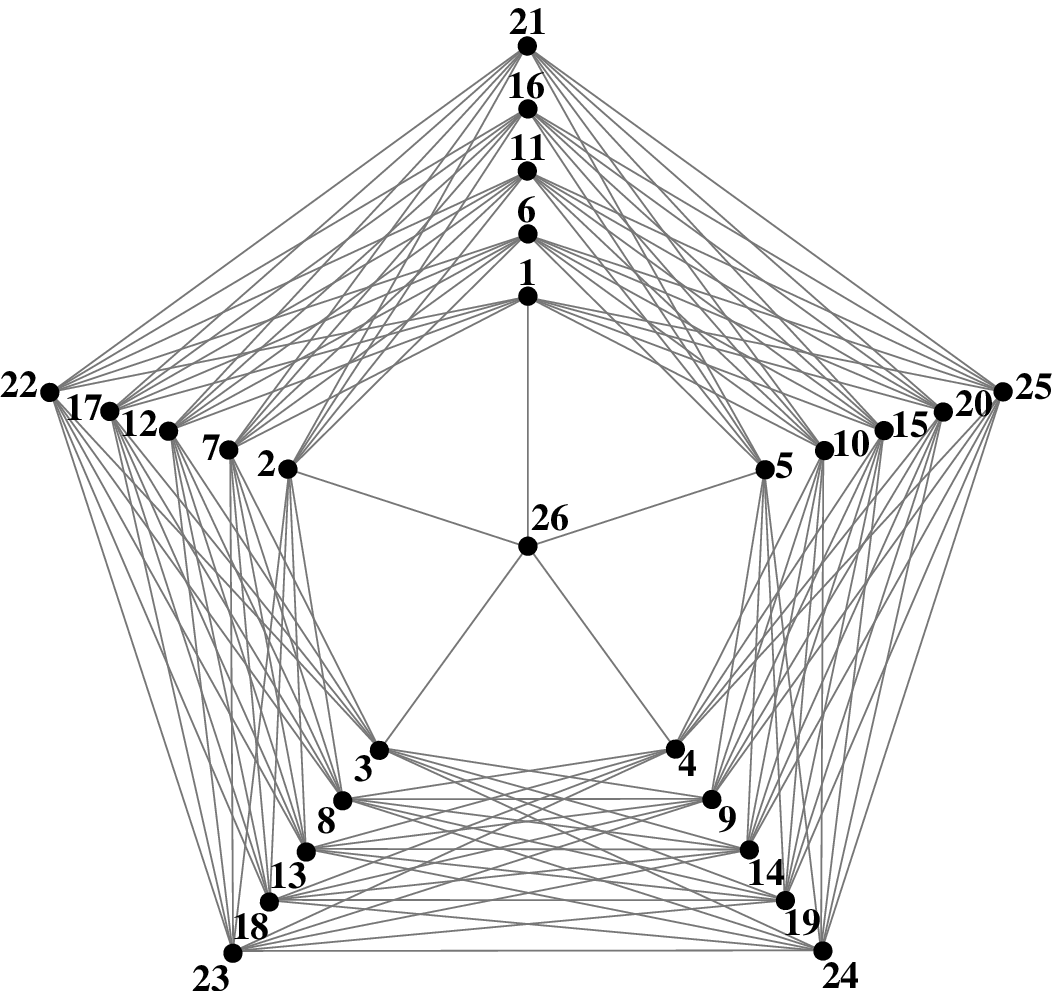} &
      \includegraphics*[width=0.393\textwidth]{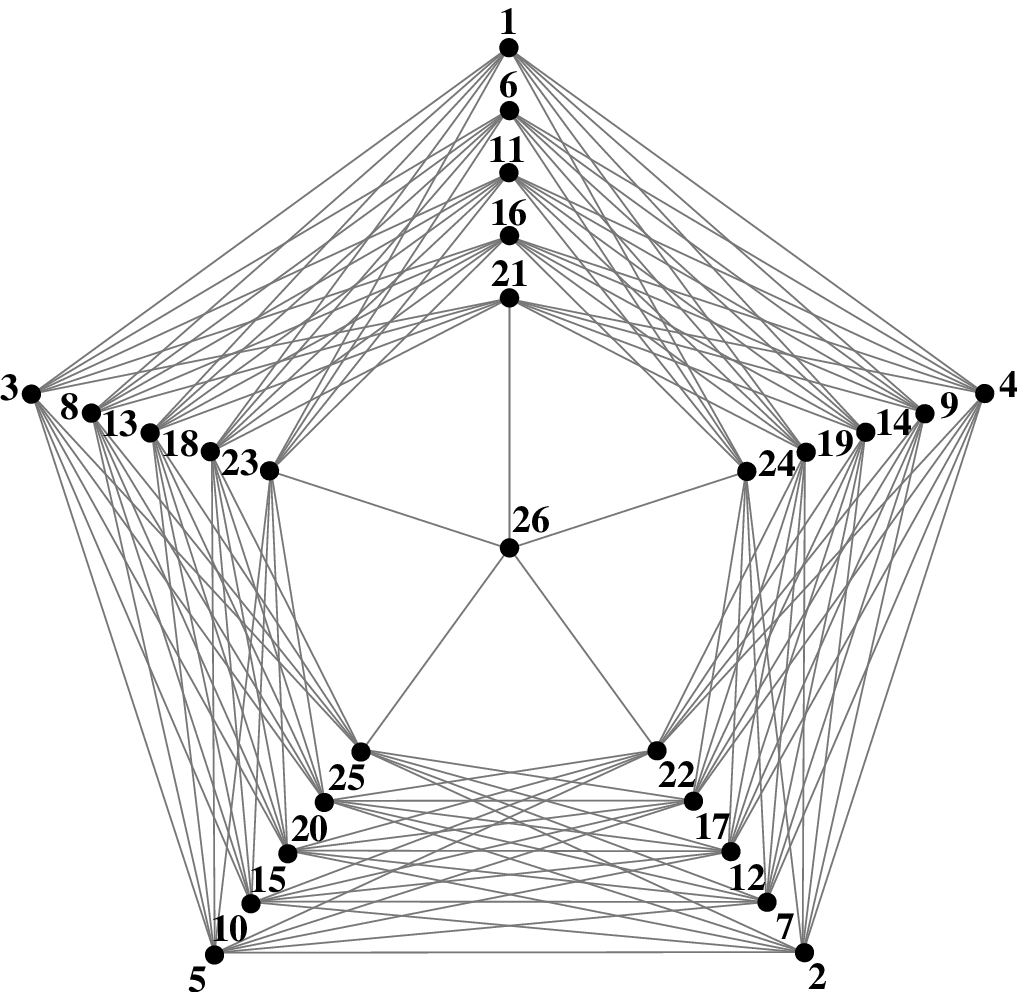}
     \end{tabular}
   \end{center}
  \caption{The graphs with adjacency matrices $M_1$ and $M_2$.}
  \label{fig:noK3}
\end{figure}

We denote the vertices from the top of the matrices as $1,2,
\ldots 26$. The graphs are shown in Figure~\ref{fig:noK3}. Vertices $1-25$ are
the triangle-free constructions from $C_5$ and $K_{5,5}-e$. Vertex
$26$ is in no triangle, since its neighbors have no edges between
them. Hence, the graphs are triangle-free.

\begin{figure}[!ht]
   \begin{center}
     \begin{tabular}{cc}
      \includegraphics*[width=0.393\textwidth]{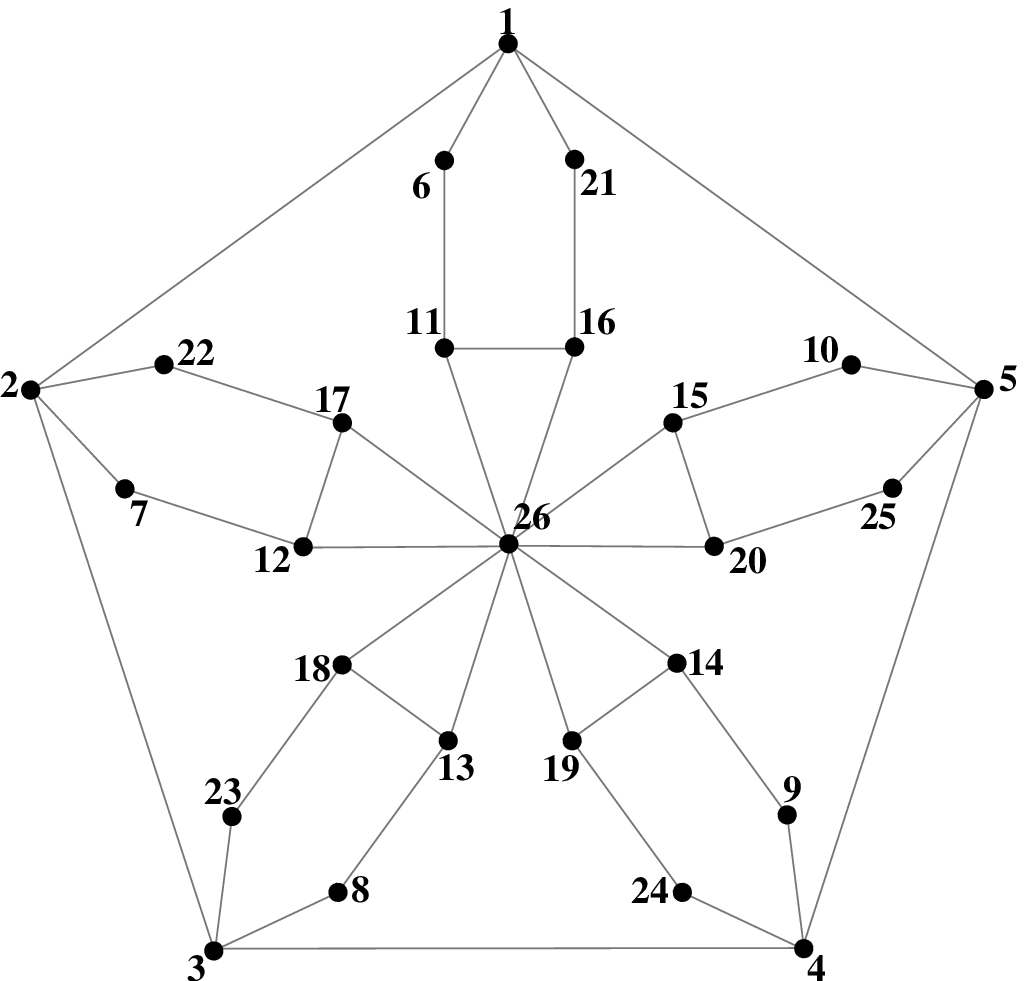} &
      \includegraphics*[width=0.407\textwidth]{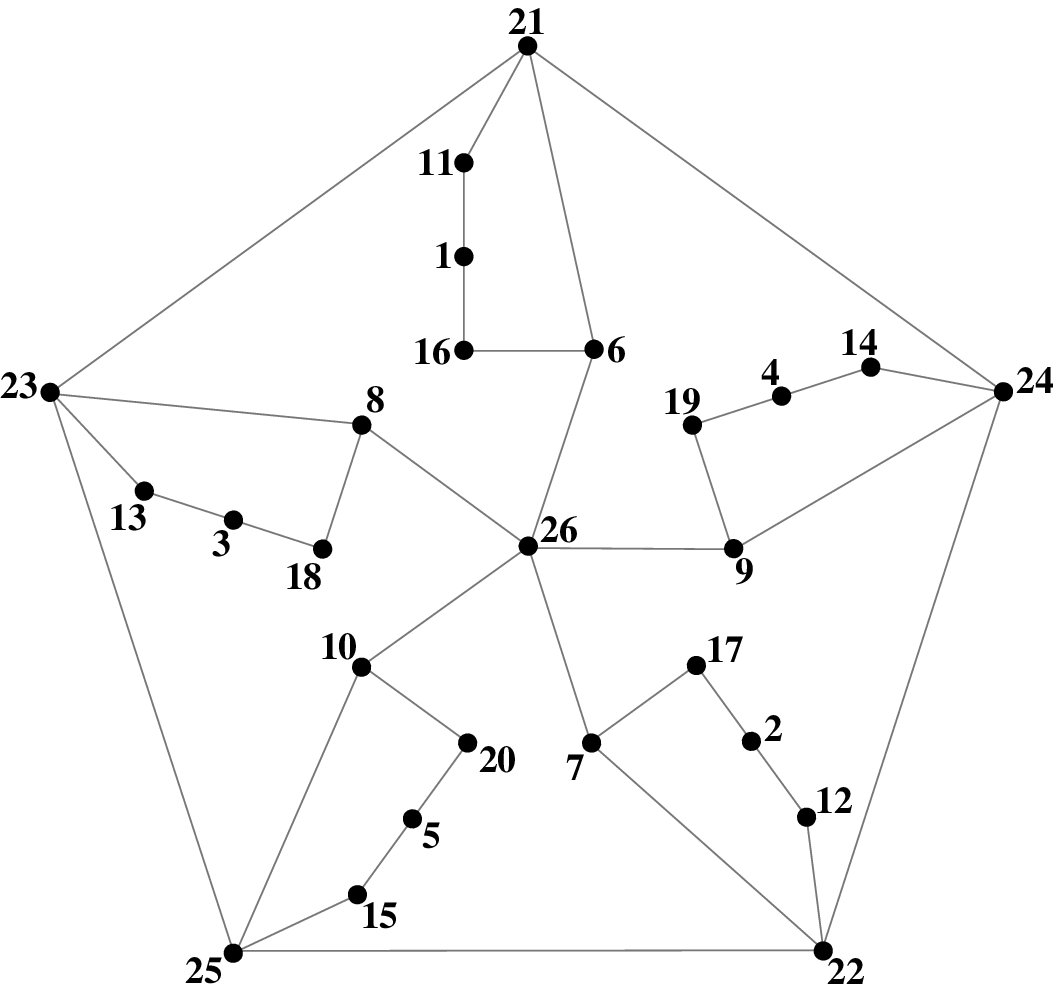}
     \end{tabular}
   \end{center}
  \caption{The graphs with adjacency matrices $M_3$ and $M_4$.}
  \label{fig:noC4}
\end{figure}

The remaining edges are distributed as described by the adjacency
matrices $M_3$ and $M_4$.
\[
\begin{array}{cc}
M_3 = \left[ \begin{array}{cccccc}
X & I & 0 & 0 & I & \bar{0} \\
I & 0 & I & 0 & 0 & \bar{0} \\
0 & I & 0 & I & 0 & \bar{1} \\
0 & 0 & I & 0 & I & \bar{1} \\
I & 0 & 0 & I & 0 & \bar{0} \\
\bar{0}^T & \bar{0}^T & \bar{1}^T & \bar{1}^T & \bar{0}^T & 0
\end{array} \right]
&
M_4 = \left[ \begin{array}{cccccc}
0 & 0 & I & I & 0 & \bar{0} \\
0 & 0 & 0 & I & I & \bar{1} \\
I & 0 & 0 & 0 & I & \bar{0} \\
I & I & 0 & 0 & 0 & \bar{0} \\
0 & I & I & 0 & Y & \bar{0} \\
\bar{0}^T & \bar{1}^T & \bar{0}^T & \bar{0}^T & \bar{0}^T & 0
\end{array} \right]
\end{array}
\]
It is not hard to see that $M_1+M_2+M_3+M_4$ is the adjacency 
matrix of $K_{26}$. It is clear from Figure~\ref{fig:noC4} that 
there are no quadrilaterals. $\square$

\newpage
\section{The Ramsey number $R_4(C_4)$}

From the Ramsey number survey \cite{r} we also know that $18\leq
R_4(C_4)\leq 21$. It was shown by Clapham, Flockhart and Sheehan
\cite{frs} that a $C_4$-free graph with $19$ vertices has at most
$42$ edges. Since $4\cdot 42=168$ and there are 171 edges in
$K_{19}$, it is not possible to four-color the edges of $K_{19}$
without a monochromatic quadrilateral.

\begin{theorem}
$R_4(C_4) \leq 19$.  
\end{theorem}

\end{document}